\documentclass[12pt]{article}%
\usepackage{sw20bams}
\usepackage{amsmath}
\usepackage{amsfonts}
\usepackage{amssymb}
\usepackage{graphicx}%
\providecommand{\U}[1]{\protect\rule{.1in}{.1in}}
\begin{document}

\title{0-Pierced Triangles within a Poisson Overlay}
\author{Steven Finch}
\date{April 24, 2019}
\maketitle

\begin{abstract}
Let the Euclidean plane be simultaneously and independently endowed with a
Poisson point process and a Poisson line process, each of unit intensity.
\ Consider a triangle $T$ whose vertices all belong to the point process.
\ The triangle is 0-pierced if no member of the line process intersects any
side of $T$. \ Our starting point is Ambartzumian's 1982\ joint density for
angles of $T$; our exposition is elementary and raises several unanswered questions.

\end{abstract}

\footnotetext{Copyright \copyright \ 2018 by Steven R. Finch. All rights
reserved.}A triangle with angles $\alpha$, $\beta$, $\gamma$ is \textbf{acute}
if $\max\{\alpha,\beta,\gamma\}<\pi/2$ and \textbf{well-conditioned} if
$\min\{\alpha,\beta,\gamma\}>\pi/6$. \ Given a random mechanism for generating
triangles in the plane, we dutifully calculate corresponding probabilities out
of sheer habit and for the sake of numerical concreteness. \ 

Beginning with a Poisson point process of unit intensity, let us form a
triangle by taking the convex hull of three particles (members of the
process). \ The triangle is \textbf{0-filled} if no other particles are
contained in the convex hull. \ Study of such configurations is complicated by
the prevalence of long, narrow triangles with angles typically $\approx0$ or
$\approx\pi$. \ We defer discussion of these until later.

Beginning with a Poisson point process and a Poisson line process, also of
unit intensity and independent, let us form a triangle as before. \ The
triangle is \textbf{0-pierced} if the intersection of each line with the
convex hull is always empty. \ Nothing is presumed about the existence or
number of other interior particles; there may be $0$ or $1$ or $2$ or many
more. \ Since the angles satisfy $\alpha+\beta+\gamma=\pi$, we can eliminate
$\gamma$ from consideration and write the joint density for $\alpha$, $\beta$
as \cite{A1-zrpctr, A2-zrpctr, K-zrpctr, O-zrpctr}%
\[
\frac{42}{\pi}\dfrac{\sin(x)\sin(y)\sin(x+y)}{\left[  \sin(x)+\sin
(y)+\sin(x+y)\right]  ^{4}}%
\]
where $x>0$, $y>0$, $x+y<\pi$. This is a remarkable result, owing to the
scattered complexity of particles overlaid with lines. Integrating out $y$, we
obtain the marginal density for $\alpha$:
\[
f(x)=\frac{21}{2\pi}\frac{\left(  7-5\cos(x)\right)  \left(  1+\cos(x)\right)
+4\left(  5-\cos(x)\right)  \left(  1-\cos(x)\right)  \ln\left(  \sin(\frac
{x}{2})\right)  }{\left(  1+\cos(x)\right)  ^{4}}%
\]
and
\[%
\begin{array}
[c]{ccc}%
\mathbb{E}(\alpha)=\dfrac{\pi}{3}=\,1.0471975511..., &  & \mathbb{E}%
(\alpha^{2})=\dfrac{13}{10}-2\ln(2)+4\ln(2)^{2}=1.8355176945....
\end{array}
\]
Corresponding to the density for $\max\{\alpha,\beta,\gamma\}$, the expression
$3f(x)$ holds when $\pi/2<x<\pi$; the expression when $\pi/3<x<\pi/2$ is
$3\tilde{f}(x)$ where%
\[
\tilde{f}(x)=\frac{21}{2\pi}\frac{\left(  7-4\cos(x)\right)  \left(
1-2\cos(x)\right)  -4\left(  5-\cos(x)\right)  \left(  1-\cos(x)\right)
\ln\left(  2\sin(\frac{x}{2})\right)  }{\left(  1+\cos(x)\right)  ^{4}}.
\]
It thus follows that%
\[
\mathbb{P}(\text{a typical 0-pierced triangle is acute})=\frac{96-132\ln
(2)-\pi}{2\pi}=0.2169249267....\text{ }%
\]
Corresponding to the density for $\min\{\alpha,\beta,\gamma\}$, the expression
$-3\tilde{f}(x)$ holds when $0<x<\pi/3$ and hence%
\[%
\begin{array}
[c]{l}%
\mathbb{P}(\text{a typical 0-pierced triangle is well-conditioned})=\medskip\\
\dfrac{-3144+1584\sqrt{3}+(-2190+1338\sqrt{3})\ln(2)+(4380-2676\sqrt{3}%
)\ln(-1+\sqrt{3})+(71+41\sqrt{3})\pi}{2(71+41\sqrt{3})\pi}\\
=0.2393922701....
\end{array}
\]
From $\mathbb{V}(\alpha+\beta+\gamma)=0$, we deduce that $\mathbb{E}%
(\alpha\beta)-\mathbb{E}(\alpha)\mathbb{E}(\beta)=-(1/2)\mathbb{V}(\alpha)$
and therefore
\[
\mathbb{E}(\alpha\,\beta)=-\dfrac{13}{10}+\ln(2)-2\ln(2)^{2}+\frac{\pi^{2}}%
{6}=0.7271752195....
\]
Simulation provides compelling evidence that Ambartzumian's \cite{A1-zrpctr,
A2-zrpctr}\ joint density is valid -- see Figure 1 -- although it does not
provide insight leading to an actual proof.%
\begin{figure}[ptb]%
\centering
\includegraphics[
height=3.0234in,
width=2.9551in
]%
{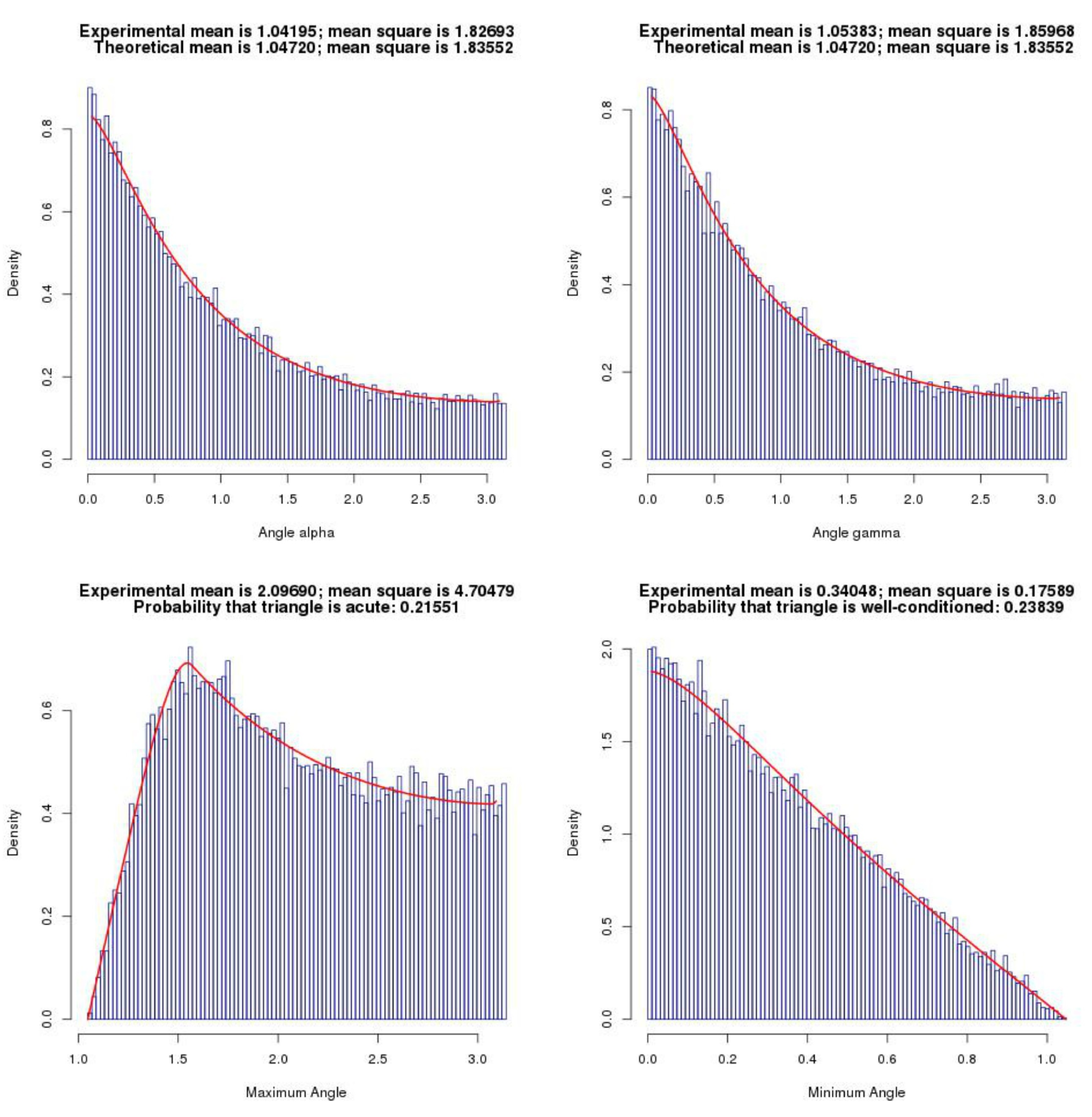}%
\caption{Histograms for angles, maximum angle and minimum angle in 0-pierced
triangles.}%
\end{figure}

\section{Related Expressions}

We turn attention to the bivariate densities%
\[
C_{j}\dfrac{\sin(x)\sin(y)\sin(x+y)}{\left[  \sin(x)+\sin(y)+\sin(x+y)\right]
^{j}}%
\]
where $x>0$, $y>0$, $x+y<\pi$ and%
\[%
\begin{array}
[c]{ccccc}%
C_{1}=\dfrac{4}{12-\pi^{2}}, &  & C_{2}=\dfrac{1}{\left(  -2+3\ln(2)\right)
\pi}, &  & C_{3}=8.
\end{array}
\]
The case $j=1$ appears in \cite{M1-zrpctr, F1-zrpctr} with regard to cells of
a Goudsmit-Miles tessellation (sampled until a triangle emerges); the case
$j=3$ appears in \cite{F2-zrpctr, F3-zrpctr} with regard to triangles created
via breaking a line segment (in two places at random). \ For $j=2$, we obtain
the univariate density for $\alpha$:
\[
g(x)=\dfrac{1}{2\left(  -2+3\ln(2)\right)  \pi}\frac{\left(  \pi-x\right)
\sin(x)+4\left(  1-\cos(x)\right)  \ln\left(  \sin(\frac{x}{2})\right)
}{1+\cos(x)}%
\]
and%
\[
\mathbb{E}(\alpha^{2})=\dfrac{4\left(  \pi^{2}-12\ln(2)\right)  \ln
(2)-3\zeta(3)}{6\left(  -2+3\ln(2)\right)  }=1.4611131303...
\]
where $\zeta(3)$ is Ap\'{e}ry's constant \cite{F4-zrpctr}.

Corresponding to the density for $\max\{\alpha,\beta,\gamma\}$, the expression
$3g(x)$ holds when $\pi/2<x<\pi$; the expression when $\pi/3<x<\pi/2$ is
$3\tilde{g}(x)$ where%
\[
\tilde{g}(x)=\dfrac{1}{2\left(  -2+3\ln(2)\right)  \pi}\frac{\left(
3x-\pi\right)  \sin(x)-4\left(  1-\cos(x)\right)  \ln\left(  2\sin(\frac{x}%
{2})\right)  }{1+\cos(x)}.
\]
It thus follows that%
\begin{align*}
\mathbb{P}(\text{a typical such triangle is acute})  &  =\frac{-24\ln
(2)+\left(  4-3\ln(2)\right)  \pi+12G}{4\left(  -2+3\ln(2)\right)  \pi}\\
&  =0.3903338870...
\end{align*}
and $G$ is Catalan's constant \cite{F5-zrpctr}. \ Corresponding to the density
for $\min\{\alpha,\beta,\gamma\}$, the expression $-3\tilde{g}(x)$ holds when
$0<x<\pi/3$ and hence%
\[%
\begin{array}
[c]{l}%
\mathbb{P}(\text{a typical such triangle is well-conditioned})=0.4190489201...
\end{array}
\]
(exact expression omitted for reasons of length). \ As before, we deduce that
\[
\mathbb{E}(\alpha\,\beta)=\dfrac{2\left(  \pi^{2}+24\ln(2)\right)  \ln
(2)-4\pi^{2}+3\zeta(3)}{12\left(  -2+3\ln(2)\right)  }=0.9143775016....
\]
What's missing, of course, is a natural procedure for generating (not
necessarily planar) triangles whose angles $\alpha$, $\beta$, $\gamma$ obey
the distributional law prescribed by $j=2$.

\section{0-Filled Triangles}

The phrase \textquotedblleft0-filled\textquotedblright\ first appeared in
\cite{M2-zrpctr, C-zrpctr}. \ Let us initially discuss the simulation
underlying 0-pierced triangles. \ Given a parameter value $\lambda>0$, we
generated data $(\alpha_{1},\beta_{1})$, $(\alpha_{2},\beta_{2})$, \ldots,
$(\alpha_{n},\beta_{n})$ via Poisson overlays in the planar disk of radius
$\lambda$ centered at the origin. Our goal was to verify a probability
theoretic expression:%
\[
\mathbb{P}(x<\alpha\leq x+dx,\;y<\beta\leq y+dy)=\frac{42}{\pi}\dfrac
{\sin(x)\sin(y)\sin(x+y)}{\left[  \sin(x)+\sin(y)+\sin(x+y)\right]  ^{4}%
}\,dx\,dy
\]
as $\lambda\rightarrow\infty$. \ This was done simply by histogramming the
data, given large enough $n$ and $\lambda$.

For 0-filled triangles, however, we face a situation where the goal is less
tangible. \ Ambartzumian's measure theoretic expression \cite{A2-zrpctr}:%
\[
\mathbb{M}(x<\alpha\leq x+dx,\;y<\beta\leq y+dy)=\dfrac{2}{\sin(x)\sin
(y)\sin(x+y)}\,dx\,dy
\]
cannot be normalized to give a probability density (that is, encompassing unit
area). \ It follows that \cite{Su-zrpctr, A3-zrpctr}%
\[
\mathbb{M}(x<\max\{\alpha,\beta,\gamma\}\leq x+dx)=\left\{
\begin{array}
[c]{lll}%
-12\csc(x)^{2}\ln\left(  2\cos(x)\right)  &  & \text{if }\pi/3\leq x<\pi/2,\\
\infty &  & \text{if }\pi/2\leq x<\pi
\end{array}
\right.
\]
and, for $0<x<\pi/3$,
\[
\mathbb{M}(x<\min\{\alpha,\beta,\gamma\}\leq x+dx)=12\csc(x)^{2}\ln\left(
2\cos(x)\right)
\]
-- see Figures 2 and 3 -- but verification is problematic.%
\begin{figure}[ptb]%
\centering
\includegraphics[
height=3.0234in,
width=3.0234in
]%
{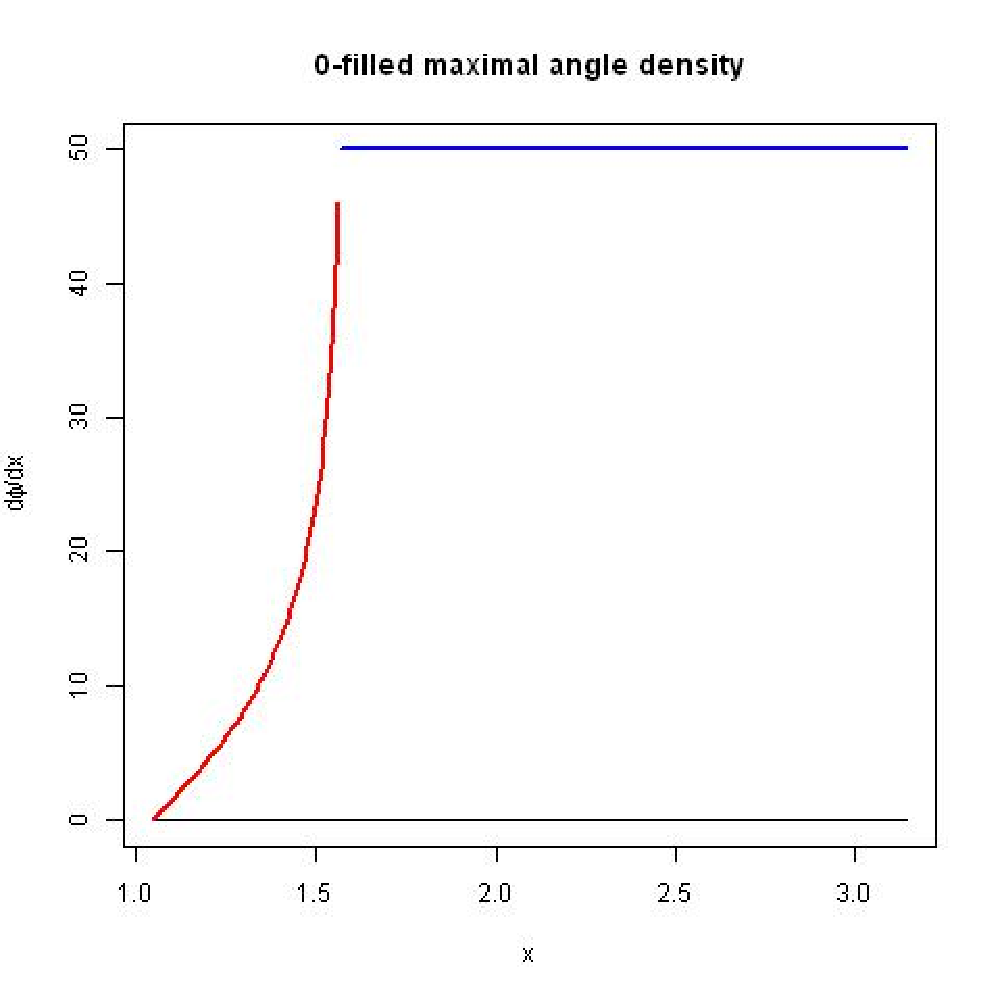}%
\caption{$\phi(x)=4\left(  3x-\pi+3\cot(x)\ln(2\cos(x))\right)  $ for
$\pi/3<x<\pi/2$.}%
\end{figure}
\begin{figure}[ptb]%
\centering
\includegraphics[
height=3.0234in,
width=3.0234in
]%
{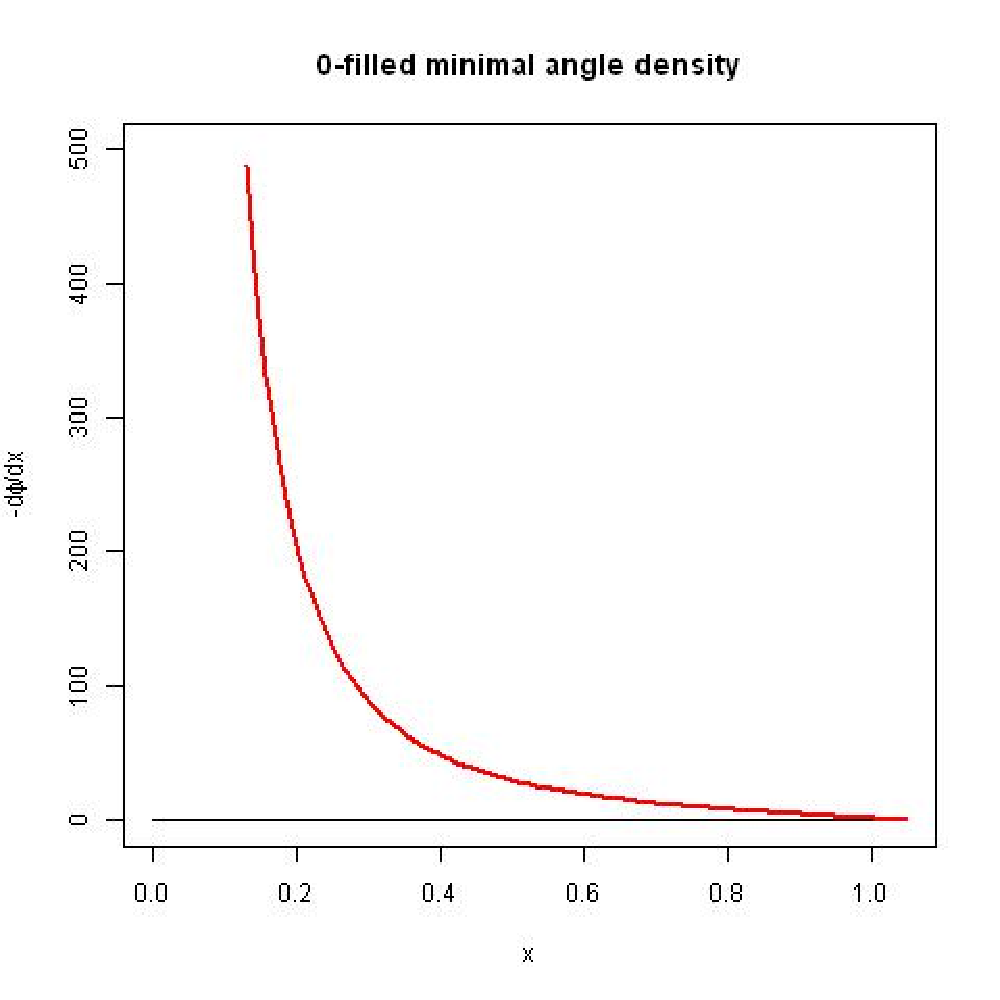}%
\caption{$\phi(x)=4\left(  3x-\pi+3\cot(x)\ln(2\cos(x))\right)  $ for
$0<x<\pi/3$.}%
\end{figure}
\ As before, we can generate data over disks of increasing radius $\lambda$.
\ Figure 4 provides histograms of $\alpha$ for $\lambda=2,3,4,5$; Figures 5
and 6 do likewise for $\max\{\alpha,\beta,\gamma\}$ and $\min\{\alpha
,\beta,\gamma\}$. \ Clearly
\[
\lim_{\lambda\rightarrow\infty}\mathbb{P}(\text{a typical 0-filled triangle is
acute})=0,
\]%
\[
\lim_{\lambda\rightarrow\infty}\mathbb{P}(\text{a typical 0-filled triangle is
well-conditioned})=0
\]
on empirical grounds. \ Unfortunately we do not know how to confirm
theoretical predictions stemming from \cite{Su-zrpctr, A3-zrpctr}:%
\[
\mathbb{M}(\text{a typical 0-filled triangle is acute})=2\pi\approx6.283,
\]%
\[
\mathbb{M}(\text{a typical 0-filled triangle is well-conditioned})=2\left(
3\sqrt{3}\ln(3)-\pi\right)  \approx5.134
\]
via our experimental simulation. \ A procedure to adjust the histogramming of
the data, in order to demonstrate an improved fit as $\lambda\rightarrow
\infty$, would be welcome.%
\begin{figure}[ptb]%
\centering
\includegraphics[
height=3.0234in,
width=2.9862in
]%
{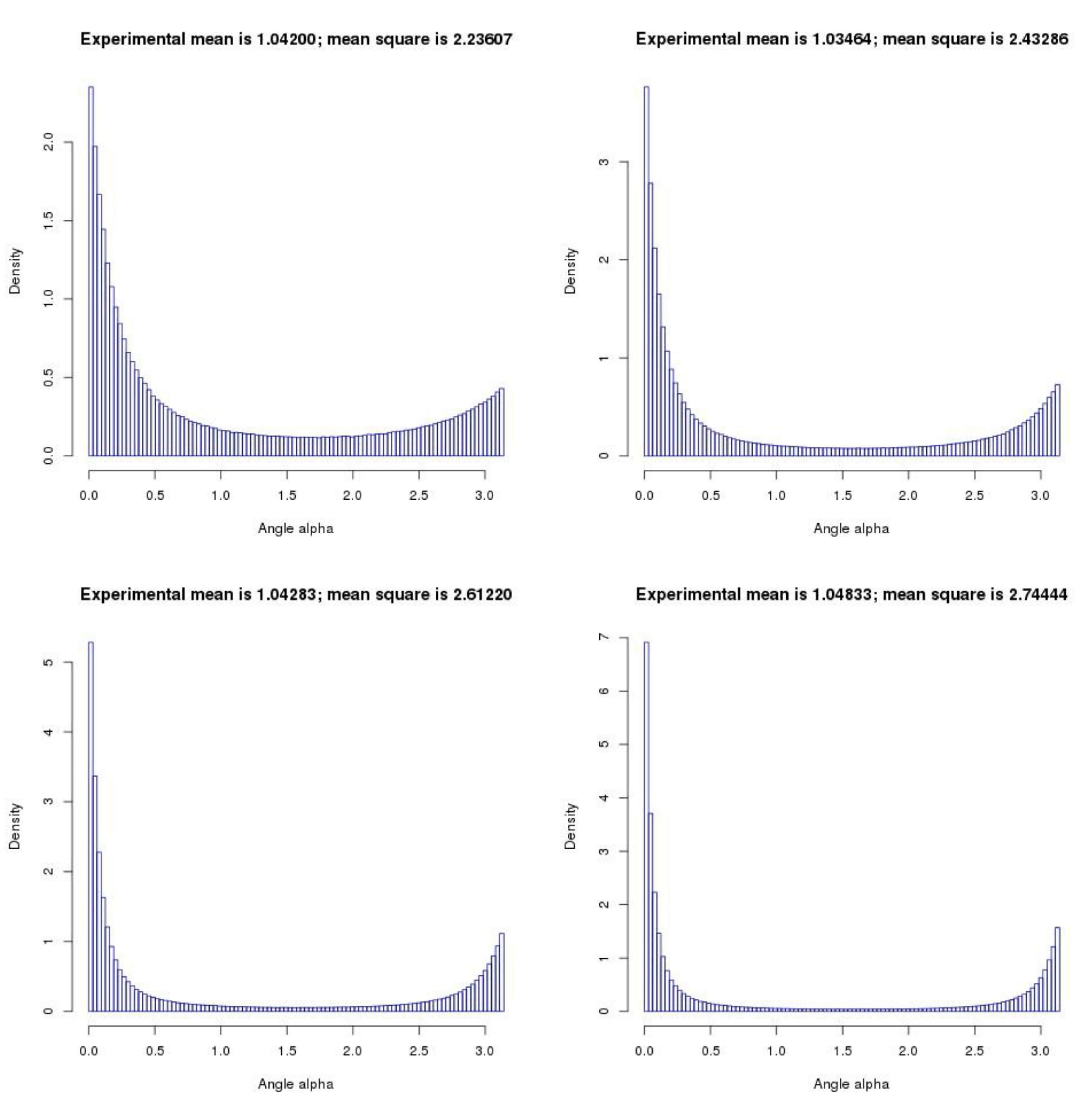}%
\caption{Histograms for arbitrary angle in 0-filled triangles, for increasing
$\lambda$.}%
\end{figure}
\begin{figure}[ptb]%
\centering
\includegraphics[
height=3.0234in,
width=2.9689in
]%
{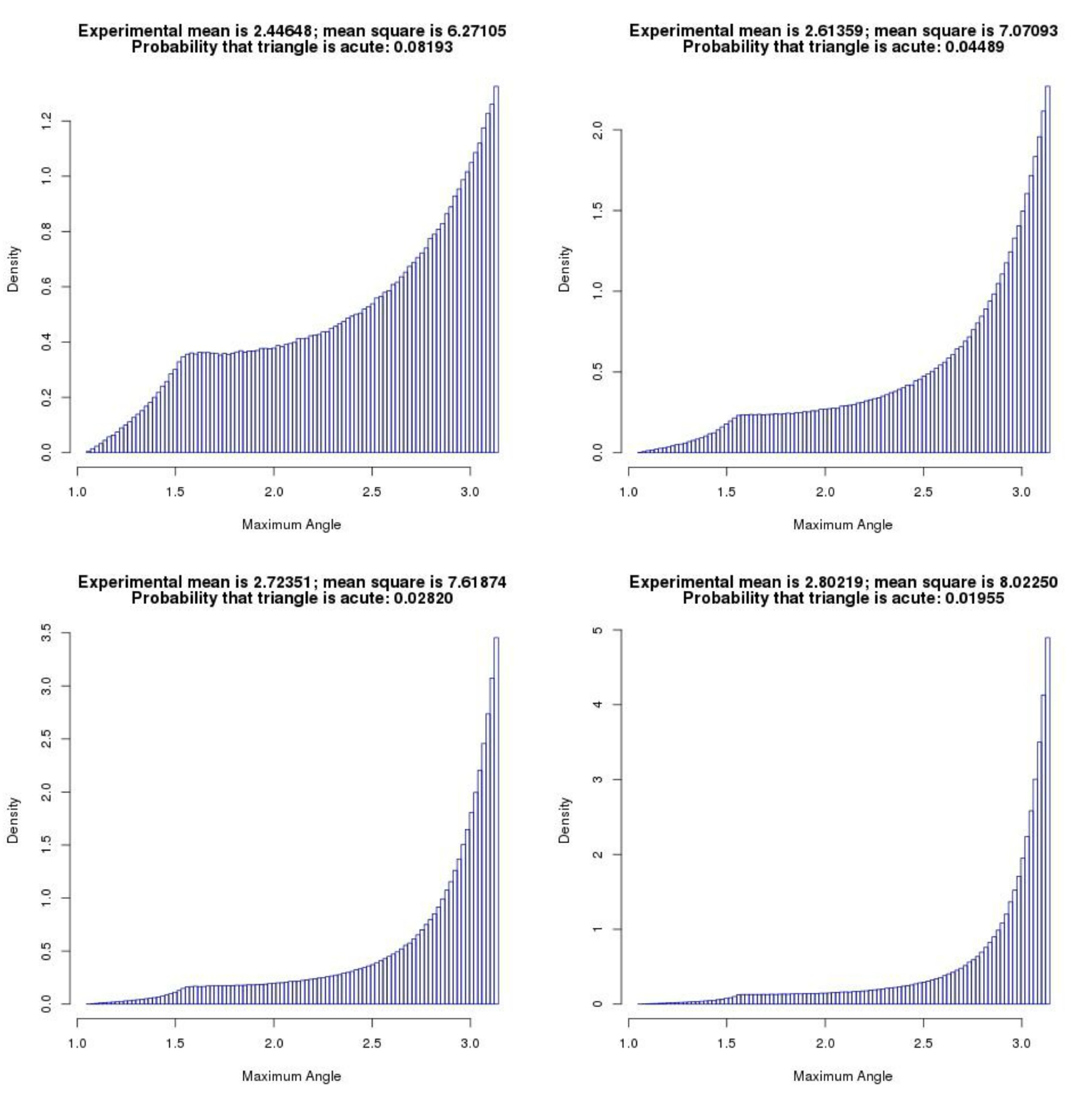}%
\caption{Histograms for maximum angle in 0-filled triangles, for increasing
$\lambda$.}%
\end{figure}
\begin{figure}[ptb]%
\centering
\includegraphics[
height=3.0234in,
width=2.9689in
]%
{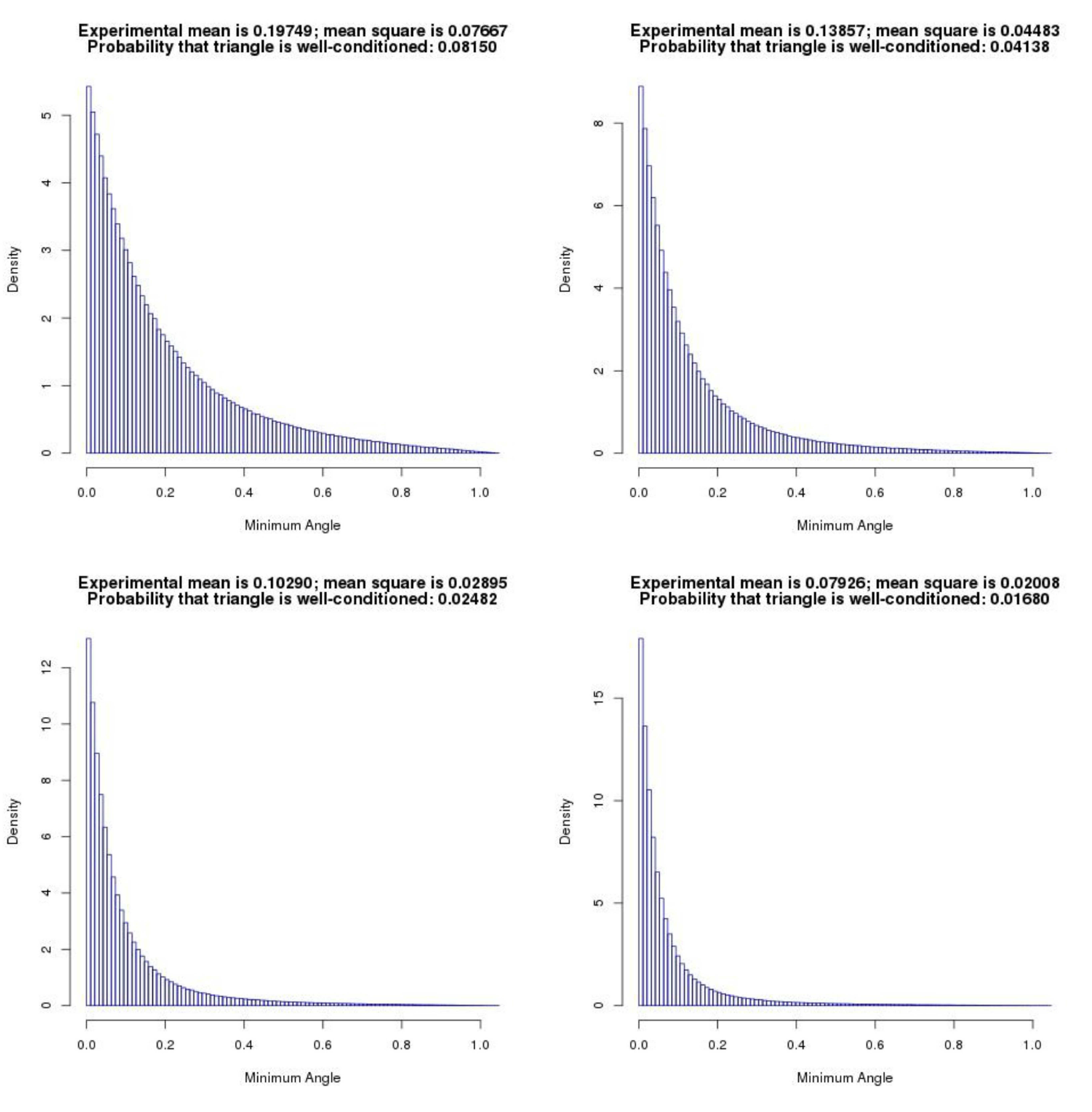}%
\caption{Histograms for minimum angle in 0-filled triangles, for increasing
$\lambda$.}%
\end{figure}

\section{Process Intensities}

We report here on work in \cite{Fa-zrpctr}. \ Given a Poisson overlay $\Omega
$, define a $\mathbf{T}_{\cdot0}$\textbf{ process} to be the set of all
0-pierced triangles within $\Omega$. \ Let the \textbf{intensity} $i$ of the
process be the mean number of triangles per unit area. \ It is known that
\[
i(T_{\cdot0})=\frac{2\pi^{2}}{21}=0.9399623239...=\frac{1}{6}%
(5.6397739434...).
\]
The factor of $1/6$ arises because the three vertices were (apparently)
ordered in \cite{Fa-zrpctr}, thus every triangle was counted $3!=6$ times.
\ We may similarly examine the set of all 0-filled triangles; it is not
surprising that $i(T_{0\cdot})=\infty$. \ Most interesting, however, is the
set of all triangles that are both 0-filled and 0-pierced:%
\begin{align*}
i(T_{00})  &  =0.6554010386...=\frac{1}{6}(3.9324062319...)\\
&  =\frac{\pi}{18\sqrt{3}}%
{\displaystyle\int\limits_{0}^{\pi}}
\frac{\xi(x)\eta(x)}{\sqrt{b(x)\left(  a(x)-c(x)\right)  }}\sin\left(
\frac{x}{2}\right)  dx
\end{align*}
where%
\[%
\begin{array}
[c]{ccccc}%
a(x)=\dfrac{2}{3}\left(  \cos\left(  \dfrac{x}{3}\right)  +1\right)  , &  &
b(x)=\dfrac{2}{3}\left(  \cos\left(  \dfrac{x-2\pi}{3}\right)  +1\right)  , &
& c(x)=\dfrac{2}{3}\left(  \cos\left(  \dfrac{x+2\pi}{3}\right)  +1\right)  ,
\end{array}
\]%
\[%
\begin{array}
[c]{ccc}%
q(x)=\sqrt{\dfrac{a(x)\left(  b(x)-c(x)\right)  }{b(x)\left(
a(x)-c(x)\right)  }}, &  & h(x)=2(27)^{1/4}\cos\left(  \dfrac{x}{2}\right)
^{-1/2},
\end{array}
\]%
\[
\xi(x)=2\left(  4+h(x)^{2}\right)  -\sqrt{\pi}\left(  6+h(x)^{2}\right)
h(x)\exp\left(  \frac{h(x)^{2}}{4}\right)  \operatorname*{erfc}\left(
\frac{h(x)}{2}\right)  ,
\]%
\[
\eta(x)=\left(  \frac{3}{c(x)}-\frac{3}{a(x)}\right)  E(q(x))+\left(  \frac
{3}{a(x)}-1\right)  K(q(x));
\]%
\[%
\begin{array}
[c]{ccc}%
K(y)=%
{\displaystyle\int\limits_{0}^{\pi/2}}
\dfrac{1}{\sqrt{1-y^{2}\sin(\theta)^{2}}}\,d\theta, &  & E(y)=%
{\displaystyle\int\limits_{0}^{\pi/2}}
\sqrt{1-y^{2}\sin(\theta)^{2}}\,d\theta
\end{array}
\]
are complete elliptic integrals of the first and second kind; and%
\[
\operatorname*{erf}(z)=\frac{2}{\sqrt{\pi}}%
{\displaystyle\int\limits_{0}^{z}}
\exp(-\tau^{2})\,d\tau=1-\operatorname*{erfc}\left(  z\right)
\]
is the error function. \ Formulas (7) and (8) in \cite{Fa-zrpctr}, devoted to
a more general scenario $T_{k\ell}$ than our $T_{00}$, specialize to%
\[%
{\displaystyle\int\limits_{0}^{\infty}}
s\exp\left(  -s-t\sqrt{s}\right)  \,ds=\frac{1}{8}\left[  2\left(
4+t^{2}\right)  -\sqrt{\pi}\left(  6+t^{2}\right)  t\exp\left(  \frac{t^{2}%
}{4}\right)  \operatorname*{erfc}\left(  \frac{t}{2}\right)  \right]
\]
for $t>0$ (avoiding use of a parabolic cylinder function $D_{-4}\left(
t/\sqrt{2}\right)  $ which is less familiar). \ 

Theory fails for $T_{00}$ -- we do not possess density predictions for the
histograms in Figure 7 -- nor do we know exact probabilities that a such a
triangle is acute or well-conditioned.%
\begin{figure}[ptb]%
\centering
\includegraphics[
height=3.0234in,
width=2.9603in
]%
{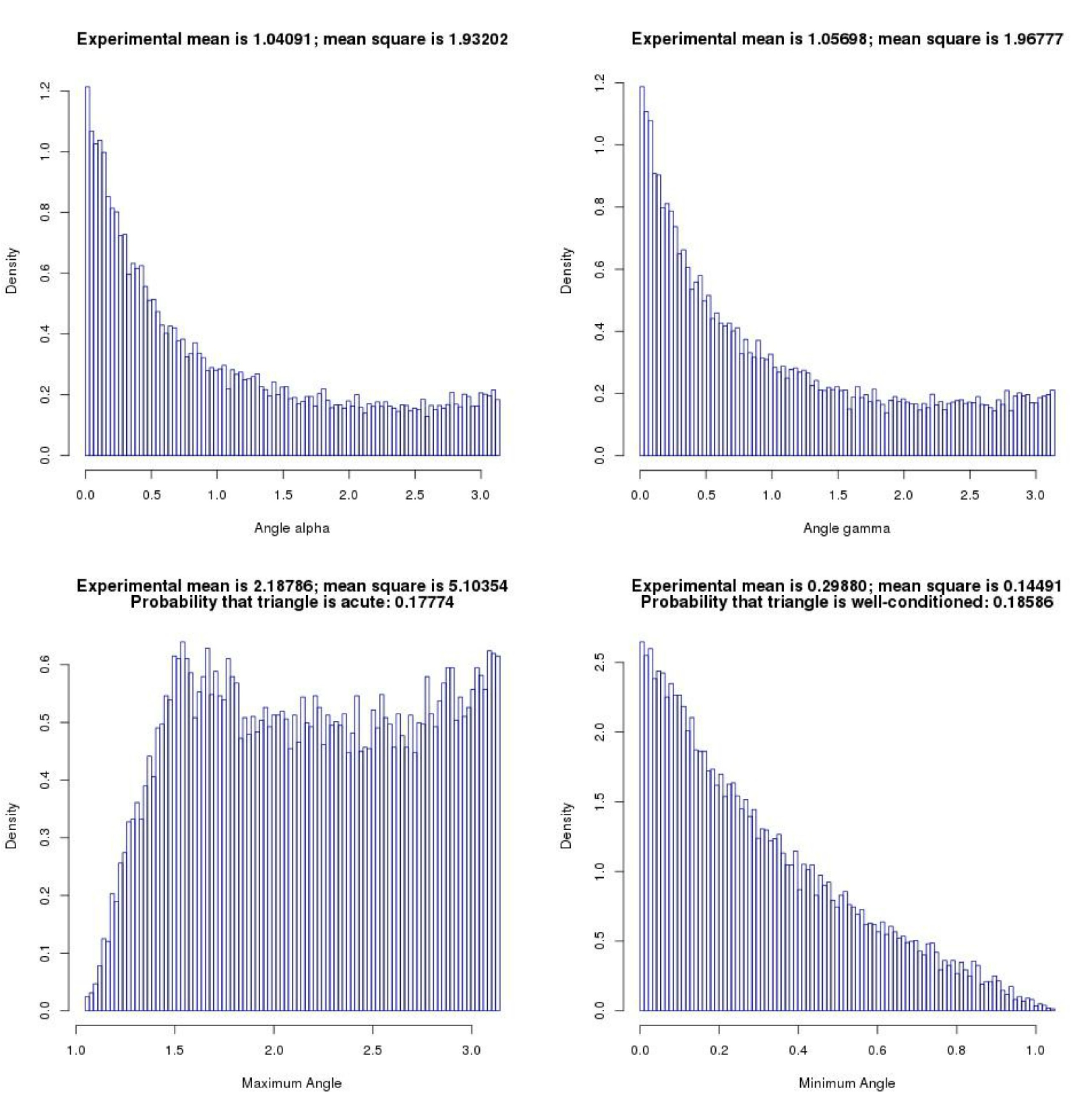}%
\caption{Histograms for angles, maximum angle and minimum angle in $T_{00}$
triangles.}%
\end{figure}


\begin{thebibliography}{99}                                                                                               %


\bibitem {A1-zrpctr}R. V. Ambartzumian, Random shapes by factorization,
\textit{Statistics in Theory and Practise, Essays in Honour of Bertil
Mat\'{e}rn}, ed. B. Ranneby, Swedish University of Agricultural Sciences,
Section of Forest Biometry, 1982, pp. 35-41; MR0688997 (84c:62004).

\bibitem {A2-zrpctr}R. V. Ambartzumian, Factorization in integral and
stochastic geometry, \textit{Stochastic Geometry, Geometric Statistics,
Stereology}, Proc. 1983 Oberwolfach conf., ed. R. Ambartzumian and W. Weil,
Teubner, 1984, pp. 14--33; MR0794864.

\bibitem {K-zrpctr}D. G. Kendall, Shape manifolds, Procrustean metrics, and
complex projective spaces, \textit{Bull. London Math. Soc.} 16 (1984) 81--121;
MR0737237 (86g:52010).

\bibitem {O-zrpctr}V. K. Oganyan, On triangle shapes formed by points of a
Poisson process in the plane (in Russian), \textit{Akad. Nauk Armyan. SSR
Dokl.}, v. 81 (1985) n. 2, 59--63; MR0826342 (87h:60032).

\bibitem {M1-zrpctr}R. E. Miles, The various aggregates of random polygons
determined by random lines in a plane, \textit{Adv. Math.} 10 (1973) 256--290;
MR0319232 (47 \#7777).

\bibitem {F1-zrpctr}S. R. Finch, Random triangles V, \textit{Mathematical
Constants II}, Cambridge Univ. Press, 2019, pp. 700--713; MR3887550.

\bibitem {F2-zrpctr}S. R. Finch, Uniform triangles with equality constraints, arXiv:1411.5216.

\bibitem {F3-zrpctr}S. R. Finch, Random triangles VI, \textit{Mathematical
Constants II}, Cambridge Univ. Press, 2019, pp. 713--718; MR3887550.

\bibitem {F4-zrpctr}S. R. Finch, Ap\'{e}ry's constant, \textit{Mathematical
Constants}, Cambridge Univ. Press, 2003, pp. 40--53; MR2003519 (2004i:00001).

\bibitem {F5-zrpctr}S. R. Finch, Catalan's constant, \textit{Mathematical
Constants}, Cambridge Univ. Press, 2003, pp. 53--59; MR2003519 (2004i:00001).

\bibitem {M2-zrpctr}R. E. Miles, On the homogeneous planar Poisson point
process, \textit{Math. Biosci.} 6 (1970) 85--127; MR0279853 (43 \#5574).

\bibitem {C-zrpctr}R. Cowan, A more comprehensive complementary theorem for
the analysis of Poisson point processes, \textit{Adv. Appl. Probab.} 38 (2006)
581--601; MR2256870 (2007k:60138).

\bibitem {Su-zrpctr}H. S. Sukiasian, Two results on triangle shapes,
\textit{Stochastic Geometry, Geometric Statistics, Stereology}, Proc. 1983
Oberwolfach conf., ed. R. Ambartzumian and W. Weil, Teubner, 1984, pp.
210--221; MR0794883.

\bibitem {A3-zrpctr}R. V. Ambartzumian, \textit{Factorization Calculus and
Geometric Probability}, Cambridge Univ. Press, 1990, pp. 61--67; 158--160;
MR1075011 (92b:60013).

\bibitem {Fa-zrpctr}V. R. Fatalov, Intensities of thinned processes of
triangles that are generated by a Poisson point process on the plane (in
Russian), \textit{Izv. Akad. Nauk Armyan. SSR Ser. Mat.}, v. 25 (1990) n. 4,
344--352, 413; Engl. transl. in \textit{Soviet J. Contemp. Math. Anal.}, v. 25
(1990) n. 4, 32--40; MR1115778 (92h:60018).

%

\begin{tabular}
[c]{lll}
& Steven Finch & \\
& MIT Sloan School of Management & \\
& Cambridge, MA, USA & \\
& \textit{steven\_finch@harvard.edu} &
\end{tabular}

\end{thebibliography}
\end{document}